\documentclass[12pt]{amsart}
\usepackage{amsmath,amsfonts,amssymb,amscd,amsthm,amsbsy,epsf}
\def\R{{\mathbb  R}}
\def\zed{{\mathbb Z}}
\begin{document}

\centerline{\bf RILEY'S CONJECTURE ON SL(2,R) REPRESENTATIONS}
\smallskip
\centerline{\bf OF 2-BRIDGE KNOTS}
\bigskip
\centerline{\sc C. McA. Gordon}

\bigskip

\subsection*{1. Introduction}	

In [R1] and [R2] Riley investigated representations of 2-bridge knot groups in $SL(2,F)$ for various fields $F$. In particular, he considered non-abelian representations in which the meridians go to parabolic elements, calling these {\it parabolic} representations. He showed that, for a given 2-bridge knot $K$, such representations correspond to the roots of a certain polynomial $\lambda_K(x) \in \zed[x]$, the {\it Riley polynomial}; see [R1, Theorem 2]. Thus the real roots of $\lambda_K(x)$ give parabolic $SL(2, \R)$ representations.
\smallskip
In [R2], Riley states ``Some of our computer calculations made in 1972-73 ... suggested that the number of real roots of [$\lambda_K(x)$] is not less than $|\sigma|/2$.'' Here $\sigma = \sigma(K)$  is the signature of $K$. Following [Tr2], we will refer to this as the 

\medskip

\noindent {\bf Riley Conjecture.}
{\it The number of real roots of the Riley polynomial of a 2-bridge knot $K$ is at least $|\sigma(K)|/2$.}

\medskip

Note that since $\lambda_K(x)$ has no multiple roots [R1, Theorem 3], the statement is unambiguous.
\smallskip

Our main result is 
\medskip

\noindent {\bf Theorem 1.1.} {\it The Riley Conjecture is true.}
\medskip

For double twist knots, the Riley Conjecture was proved by Tran [Tr2].
\smallskip

One of our interests in the Riley Conjecture is its connection with the question of when the $n$-fold cyclic branched cover $\Sigma_n(K)$ of a knot $K$ has left-orderable fundamental group. More precisely, as pointed out in [Tr2], Hu's argument in [H] shows that Theorem 1.1 has the following corollary.
\medskip

\noindent {\bf Corollary 1.2.} {\it Let $K$ be a 2-bridge knot with $\sigma(K) \ne 0$. Then $\Sigma_n(K)$ has left-orderable fundamental group for $n$ sufficiently large.}
\medskip

By contrast, there are 2-bridge knots $K$ such that $\Sigma_n(K)$ has non-left-orderable fundamental group for all $n$, by [Te, Proof of Theorem 2] and [BGW, Theorem 4].

For any knot $K$, the determinant and signature are related by the following congruence [M, Theorem 5.6]
\medskip

\hspace{130pt}$\det(K) \equiv (-1)^{\sigma(K)/2}$ (mod 4)
\medskip

If $K$ is the 2-bridge knot corresponding to $p/q \in \mathbb Q$, $p > 0$, then $\det(K) = p$. Hence if $p \equiv -1$ (mod 4) then $\sigma(K) \equiv 2$ (mod 4), and Corollary 1.2 applies. In this case the conclusion of Corollary 1.2 was proved by Hu [H]. 
 
For other results on the left-orderability of the fundamental groups of cyclic branched covers of knots see [GL] and [Tr1].
\medskip

\noindent {\bf Question 1.3.} {\it Does Corollary 1.2 hold without the assumption that $K$ is 2-bridge?}
\medskip

The proof of Theorem 1.1 uses a variant of the classical theorem of Sturm on the number of real roots of a polynomial with real coefficients. This is treated in Section 2. In Section 3 we prove the Riley Conjecture, and in Section 4 we discuss Corollary 1.2.
\smallskip

\noindent {\it Acknowledgements}. I would like to thank Steve Boyer and Anh Tran for helpful conversations. This research was partially supported by NSF Grant DMS-1309021.

\subsection*{2. Sturm's theorem}	

Sturm's theorem gives a way of determining the number of real roots of a polynomial with real coefficients; for a nice discussion of this, including some history, see [GR]. The method depends on constructing a sequence of polynomials $f_0, f_1,...,f_n = f$ with certain properties (we find it convenient to reverse the usual numbering convention). In Theorem 2.1 we prove a version of Sturm's theorem where the key properties of $f_0, f_1,...,f_{n-1}$ are as in the classical setting, but the hypothesis on the relation between $f_n$ and $f_{n-1}$ is weakened. The conclusion is then an inequality rather than an equality.

Let $\boldsymbol\alpha$ = $(\alpha_k) = (\alpha_0, \alpha_1,..., \alpha_n)$ be a sequence of non-zero real numbers. Define the {\it variation} var$(\boldsymbol \alpha)$ of $\boldsymbol\alpha$ to be the number of changes in the corresponding sequence of signs sign$(\boldsymbol\alpha) =$ (sign$(\alpha_k))$, i.e.
\begin{equation*}
\textup{var}(\boldsymbol \alpha) = \#\{k : \alpha_k \alpha_{k+1} < 0,\: 0 \le k < n\} 
\end{equation*}

Let $\boldsymbol{f} = (f_0,f_1,...,f_n)$ be a sequence of polynomials in $\mathbb{R}[X]$. If $x \in \mathbb{R}$, set $\boldsymbol{f}(x) = (f_0(x),f_1(x),...,f_n(x)) \in \mathbb{R}^{n+1}$.

Let $Z_k = \{$real roots of $f_k \} \subset \mathbb{R}, \: 0 \le k \le n$, and let $Z = \bigcup\limits_{k=0}^n Z_k$. Choose $x_+$ and $x_{-} \in \mathbb{R}$ such that $Z \subset (x_{-}, x_{+})$. Then sign$(\boldsymbol{f}(x_{+}))$ is independent of the choice of $x_+$, so we denote it by sign$(\boldsymbol{f}(\infty))$, and write var$(\boldsymbol{f}(\infty)) =$ var$(\boldsymbol{f}(x_{+}))$. Similarly, we write sign$(\boldsymbol{f}(-\infty)) =$ sign$(\boldsymbol{f}(x_{-}))$ and var$(\boldsymbol{f}(-\infty)) =$ var$(\boldsymbol{f}(x_{-}))$.


\bigskip

\noindent{\bf Theorem 2.1.} {\it Let $\boldsymbol{f} = (f_0,f_1,...,f_n)$ be a sequence of polynomials in $\mathbb{R}(X)$ such that 

(1) $f_0$ is constant and non-zero, and 

(2) if $f_k(x_0) = 0$ for some $0<k<n$ and $x_0 \in \mathbb{R}$, then $f_{k-1}(x_0)f_{k+1}(x_0) < 0$.

\noindent Then $f_n$ has at least $|\textup{var}(\boldsymbol{f}(-\infty)) - \textup{var}(\boldsymbol{f}(\infty))|$ distinct real roots.} 

\begin{proof}
The theorem is trivially true when $n = 0$ so we assume $n \ge 1$.
\smallskip

Define $v: \mathbb{R}\setminus Z \to \{0,1,...,n\}$ by $v(x) =$ var$(\boldsymbol{f}(x))$. Note that $v$ is constant on each component of $\mathbb{R} \setminus Z$.

Suppose $x_0 \in Z$, so $x_0 \in Z_k$ for some $k$ (not necessarily unique) with $1\le k\le n$.

If $k<n$ then by condition (2) there exists $\delta > 0$ such that $f_{k-1}(x)f_{k+1}(x) < 0$ for all $x \in (x_0-\delta,x_0+\delta)$. Hence, as $x$ passes through $x_0$ the signs of $(f_{k-1}(x),f_k(x),f_{k+1}(x))$ change as $(\pm, \epsilon, \mp) \to (\pm, \epsilon',\mp)$, where $\epsilon, \epsilon' \in \{+,-\}$ This contributes 0 to the change in $v(x)$.

Suppose $k=n$. Note that $f_{n-1}(x_0) \ne 0$, by (1) if $n=1$ and by (2) if $n >1$. Hence, as $x$ passes through $x_0$ the signs of $(f_{n-1}(x),f_n(x))$ change as $(\pm,\epsilon) \to (\pm, \epsilon')$. Thus the corresponding change in $v(x)$ is 0 or $\pm1$.

Therefore $|\textup{var}(\boldsymbol{f}(-\infty)) - \textup{var}(\boldsymbol{f}(\infty))|$ is at most the number of distinct real roots of $f_n$.
\end{proof}

\subsection*{3. The Riley Conjecture}
 
Let $K$ be the 2-bridge knot corresponding to $p/q \in \mathbb{Q}$, where $p$ and $q$ are coprime, and $p$ is odd and $>1$. Let $n = (p-1)/2$. Then (see [R1, Proposition 1]) there exist $\epsilon_i, \eta_i \in \{\pm1\}$, with $\epsilon_i = \eta_{n+1-i}, \: 1 \le i \le n$, such that $\pi(K) = \pi_1(S^3 \setminus K)$ has presentation
\medskip

\hspace{120pt}$<a,b : wa = bw>$, 
\medskip

where $a$ and $b$ are meridians and $w = \prod\limits_{i = 1}^{n}a^{\epsilon_i}b^{\eta_i}$.
\smallskip

Also, $\sigma(K) = \sum\limits_{i=1}^{n}(\epsilon_i + \eta_i)$ [S]. Hence $\sum\limits_{i=1}^{n}\epsilon_i = \sigma(K)/2$.
\medskip

Riley considers parabolic representations of $\pi(K)$ into $SL(2,\mathbb{C})$, where
\medskip

$a \to \begin{pmatrix}1&1\\0&1\end{pmatrix} = A$, and
\smallskip

$b \to \begin{pmatrix}1&0\\x&1\end{pmatrix} = X$.
\medskip

Let $W_k = \prod\limits_{i=1}^{k}A^{\epsilon_{i}}X^{\eta_{i}}, \: 1 \le k \le n$, and set $W_0 = I$.

Write $W_k = \begin{pmatrix}a_k & b_k\\\ast & \ast\end{pmatrix}, \: a_k, b_k \in \mathbb {Z}[x]$, $\: 0 \le k \le n$.
\smallskip

The {\it Riley polynomial} of $K$ is defined to be $\lambda_K = a_n$. Riley showed [R1, Theorem 2] that the above assignment of $a$ and $b$ defines a homomorphism from $\pi(K)$ to $SL(2,\mathbb C)$ if and only if $\lambda_{K}(x) = 0$. Thus the real roots of $\lambda_K$ give parabolic representations of $\pi(K)$ into $SL(2, \mathbb R)$. 
\smallskip

Let $\delta_i = \epsilon_i \eta_i, \: 1\le i \le n$. Then
\medskip

$A^{\epsilon_i}X^{\eta_i} = \begin{pmatrix}{1+\delta_{i}x}&\epsilon_i\\{\eta_i}&1\end{pmatrix}$,
\smallskip

\noindent giving the recurrence equations, for $1 \le k \le n$,
\begin{equation*}
a_k = (1+\delta_{k}x)a_{k-1} + (\eta_k x)b_{k-1}
\tag{3.1}
\end{equation*}
\begin{equation*}
b_k = \epsilon_{k}a_{k-1} + b_{k-1}
\tag{3.2}
\end{equation*}
\smallskip
It follows from (3.1) and (3.2) by induction on $k$ that $a_k$ has degree $k$, with leading coefficient $\prod\limits_{i=1}^{k}\delta_i$, and 

\begin{equation*}
a_k(0) = 1
\tag{3.3}
\end{equation*}

\noindent Also, since det$W_k = 1$, we have that for all $x \in \mathbb C$, 
\bigskip

\noindent (3.4)\hspace{100pt}$a_k (x)$ and $b_k (x)$ are not both zero, $0 \le k \le n$

\bigskip

\noindent{\bf Lemma 3.1.}
{\it If $0 < k < n$ and $a_k (x_0) = 0$, $x_0 \in \mathbb R$, then $a_{k-1}(x_0)$ and $a_{k+1}(x_0)$ are non-zero and \textup{sign}$(a_{k-1}(x_0))$\textup{sign}$(a_{k+1}(x_0)) = - \eta_k \eta_{k+1}$.}

\begin{proof}
Suppose $a_k(x_0) = 0$. Then (3.1) gives
\begin{equation*}
(1+\delta_{k} x_0)a_{k-1}(x_0) + (\eta_k x_0)b_{k-1}(x_0) = 0
\tag{3.5}
\end{equation*}
\medskip
while from (3.2) we get 
\begin{equation*}
b_k (x_0) = \epsilon_k a_{k-1}(x_0) + b_{k-1}(x_0)
\tag{3.6}
\end{equation*}
\medskip
\noindent Multiplying both sides of (3.6) by $\eta_k x_0$ and using (3.5) gives
\begin{equation*}
a_{k-1}(x_0) = - (\eta_k x_0)b_k (x_0)
\tag{3.7}
\end{equation*}
\medskip
Replacing $k$ by $k+1$ in (3.1) we obtain
\begin{equation*}
a_{k+1}(x_0) = (\eta_{k+1} x_0)b_k (x_0)
\tag{3.8}
\end{equation*}
\medskip
By (3.3) $x_0 \ne 0$, and by (3.4) $b_k (x_0) \ne 0$.
The result now follows from (3.7) and (3.8).
\end{proof}
\smallskip

\begin{proof}[Proof of Theorem 1.1.] Define $f_k = (\prod\limits_{i=1}^k\eta_i)a_k,
\: 0\le k \le n$. Then $f_0$ is the constant polynomial $1$, and Lemma 3.1 implies that if $f_k (x_0) = 0$ for some $0 < k < n$ then $f_{k-1}(x_0)f_{k+1}(x_0) < 0$. Thus $\boldsymbol{f} = (f_k)$ satisfies the hypotheses of Theorem 2.1.

The coefficient of $x^k$ in $f_k$ is $(\prod\limits_{i=1}^{k}\eta_i)(\prod\limits_{i=1}^{k}\delta_i) = \prod\limits_{i=1}^{k}\epsilon_i = \mu_k \:$, say.

Then sign$(\boldsymbol {f}(\infty)) = (\mu_k)$, and sign$(\boldsymbol{f}(-\infty)) = ((-1)^{k}\mu_k)$. Since $\mu_k = \epsilon_{k}\mu_{k-1}$, we have 

\begin{equation*}
\textup{var}(\boldsymbol{f}(\infty)) = \#\{k : 1 \le k \le n, \; \epsilon_k = -1\} 
\end{equation*}

and
\begin{equation*}
\textup{var}(\boldsymbol{f}(-\infty)) = \#\{k : 1 \le k \le n, \; \epsilon_k = +1\}
\end{equation*}
\smallskip

Therefore by Theorem 2.1 the number of real roots of $\lambda_K = a_n = \pm f_n$
is at least 
\medskip
\begin{equation*}
|\textup{var}(\boldsymbol{f}(-\infty)) - \textup{var}(\boldsymbol{f}(\infty))| = |\sum\limits_{k=1}^n \epsilon_k | = |\sigma(K)|/2. 
\end{equation*}
\end{proof}
\smallskip

\noindent {\it Remark.} The inequality in the Riley Conjecture can be strict. For example, the knot $10_{32}$, which is the 2-bridge knot corresponding to the rational number $69/29$, has $\sigma(10_{32}) = 0$. On the other hand, by [KT1] and [KT2] (see also [ORS]), there is a meridian-preserving epimorphism from $\pi(10_{32})$ to $\pi(3_1)$, the group of the trefoil. Since $\pi(3_1)$ has a real parabolic representation [R1], so does $10_{32}$. 

This also shows that the converse of Corollary 1.2 is not true. In fact, by [GL, Theorem 1.2 and Lemma 9.1], $\Sigma_{n}(10_{32})$ has left-orderable fundamental group for $n \ge 6$.

\subsection*{4. Cyclic branched covers}

In this section we indicate how the argument in [H] gives Corollary 1.2.

In [R2] Riley considers arbitrary non-abelian $SL(2, \mathbb C)$ representations of $\pi(K)$, $K$ a 2-bridge knot. Up to conjugation, we may assume that 
\medskip

$a \to \begin{pmatrix}t&1\\0&t^{-1}\end{pmatrix}$, and
\smallskip

$b \to \begin{pmatrix}t&0\\x&t^{-1}\end{pmatrix}$.
\medskip

Riley shows that this defines a representation if and only if $\phi(t,x) = 0$ for a certain polynomial $\phi \in \mathbb{Z}[t^{\pm1}, x]$. He notes that $\phi(t,x) = \phi(t^{-1},x)$ [R2, Proposition 1], and therefore $\phi(t,x) = \psi(s,x)$, where $s = t + t^{-1}$, for some $\psi \in \mathbb{Z}[s,x]$. Then $\psi(2,x)$ is the Riley polynomial $\lambda_K(x)$.

Suppose $\lambda_K$ has a real root $x_0$. Since $\lambda_K$ has no repeated factors [R1, Theorem 3], $\frac{\partial \psi}{\partial x}\big|_{(s=2,\, x=x_0)}$ is non-zero. It follows that there exists $\delta > 0$ and a continuous function $\gamma : (2-\delta,2+\delta) \to \mathbb R$, with $\gamma(2) = x_0$, such that $\psi(s, \gamma(s)) = 0$ for all $s \in (2-\delta,2+\delta)$. In particular, for all $s \in (2-\delta,2)$ there is a non-abelian representation $\rho_s : \pi(K) \to
SL(2, \mathbb R)$ such that $\rho_{s}(a)$ has trace $s$. Conjugating $\rho_s$ we may assume that 

\begin{equation*}
\rho_s(a) = \begin{pmatrix}\cos\,\theta&-\sin\,\theta\\\sin\,\theta&\cos\,\theta\end{pmatrix}
\end{equation*}
where $s = 2\cos\,\theta$.

For $n$ sufficiently large, $s_n = 2\,$cos$(2\pi/n) \in (2-\delta, 2)$. Then $\rho_{s_n} (a)$ has order $n$. It follows from [H, Theorem 3.1] (see also [BGW, Theorem 6]) that $\pi_1(\Sigma_{n}(K))$ is left-orderable.

\bigskip
\bigskip

{\footnotesize\parindent=0pt
Department of Mathematics \par
The University of Texas at Austin \par
1 University Station C1200 \par
Austin TX 78712-0257  USA\par
\medskip

gordon@math.utexas.edu}

\end{document}